\documentclass[a4paper,11pt]{article}
\usepackage{latexsym}
\usepackage{amssymb}
\usepackage{amsfonts}
\usepackage{amsmath}
\usepackage[dvips]{graphicx}
\usepackage{epsf}
\newcommand{\qed}{$\Box$}

\newenvironment{@abssec}[1]{%
    \if@twocolumn

      \section*{#1}%
    \else

      \vspace{.05in}\footnotesize
      \parindent .2in
 {\upshape\bfseries #1. }\ignorespaces
    \fi}

    {\if@twocolumn\else\par\vspace{.1in}\fi}

\newenvironment{keywords}{\begin{@abssec}{\keywordsname}}{\end{@abssec}}

\newenvironment{AMS}{\begin{@abssec}{\AMSname}}{\end{@abssec}}

\newcommand\keywordsname{Key words}
\newcommand\AMSname{AMS subject classifications}
\newcommand\AMname{AMS subject classification}
\newtheorem{theorem}{Theorem}
 \newtheorem{lemma}[theorem]{Lemma}

\newtheorem{remark}[theorem]{Remark}
 
\def\qed{\vbox{\hrule height0.6pt\hbox{%
  \vrule height1.3ex width0.6pt\hskip0.8ex
  \vrule width0.6pt}\hrule height0.6pt
 }}

\title{Stationary level surfaces and Liouville-type theorems characterizing hyperplanes
\thanks{This research was partially supported by a Grant-in-Aid
for Scientific Research (B) ($\sharp$ 20340031) of
Japan Society for the Promotion of Science.}}

\author{Shigeru Sakaguchi\thanks{Department of Applied Mathematics,
Graduate School of  Engineering, Hiroshima
University, Higashi-Hiroshima, 739-8527,  Japan.
({\tt sakaguch@amath.hiroshima-u.ac.jp}).}}

\date{}

\begin{document}

\maketitle

\begin{abstract}
We consider an entire graph $S: x_{N+1} = f(x), x \in \mathbb R^N$ in $\mathbb R^{N+1}$ of a continuous real function $f$ over $\mathbb R^{N}$ with $N\ge 1$. Let $\Omega$ be an unbounded domain in $\mathbb R^{N+1}$ with boundary $\partial\Omega=S$.  Consider  nonlinear diffusion equations of the form  $\partial_t U= \Delta \phi(U)$ containing the heat equation $\partial_t U= \Delta U$. Let $U=U(X,t) =U(x, x_{N+1}, t)$ be the solution of either the initial-boundary value problem over $\Omega$ where the initial value equals zero and the boundary value equals $1$,  or 
the Cauchy problem where the initial data is  the characteristic function of the set $\mathbb R^{N+1}\setminus \Omega$. The problem we consider is to characterize $S$ in such a way that there exists a stationary level surface of $U$ in $\Omega$.

 We introduce a new class $\mathcal A$ of entire graphs $S$ and, by using the sliding method due to Berestycki, Caffarelli, and Nirenberg, we show that $S\in\mathcal A$ must be a hyperplane if there exists a stationary level surface of $U$ in $\Omega$. This is an improvement of the previous result. Next, we consider the heat equation in particular and we introduce the class $\mathcal B$ of entire graphs $S$ of functions $f$ such that each $\{ |f(x)-f(y)|: |x-y| \le 1 \}$ is bounded. With the help of the theory of viscosity solutions, we show that $S \in \mathcal B$ must be a hyperplane if there exists a stationary isothermic surface of $U$ in $\Omega$. 
 This is a considerable improvement of the previous result.
 
 Related to the problem, we consider a class $\mathcal W$ of Weingarten hypersurfaces in $\mathbb R^{N+1}$ with $N \ge 1$.  Then we show that, if $S$ belongs to $\mathcal W$ in the viscosity sense and $S$ satisfies some natural geometric condition, then $S \in \mathcal B$  must be a hyperplane.
This is also a considerable improvement of the previous result.

\end{abstract}




\pagestyle{plain}
\thispagestyle{plain}
\markboth{S. SAKAGUCHI}{Stationary level surfaces and Liouville-type theorems}

\section{Introduction}
 For $f \in C(\mathbb R^N)$ where $N \ge 1$,  let $\Omega$ be a domain in $\mathbb R^{N+1}$ given by
 \begin{equation}
 \label{upper-half domain}
 \Omega = \{ X=(x, x_{N+1}) \in \mathbb R^{N+1}\ :\ x_{N+1} > f(x) \}.
 \end{equation}
 Throughout this paper we write $X = (x, x_{N+1}) \in \mathbb R^{N+1}$ for $x =(x_1, \dots, x_N) \in \mathbb R^N$. Then we notice that $\partial\Omega =\partial\left(\mathbb R^{N+1}\setminus\overline{\Omega}\right)$.
 Let $\phi : \mathbb R \to \mathbb R$ satisfy
\begin{equation}
\label{nonlinearity}
\phi \in C^2(\mathbb R), \quad \phi(0) = 0, \ \mbox{ and }\ 0 < \delta_1 \le \phi^\prime(s) \le \delta_2\ \mbox{ for } s \in \mathbb R,
\end{equation}
 where $\delta_1, \delta_2$ are positive constants. Consider the unique bounded solution $U = U(X, t)$ of either
  the  initial-boundary value problem: 
\begin{eqnarray}
&\partial_t U=\Delta \phi(U)\ \ &\mbox{in }\ \Omega\times (0,+\infty),\label{diffusion}\\
&U=1\ \ &\mbox{on }\ \partial\Omega\times (0,+\infty),\label{dirichlet}\\
&U=0\ \ &\mbox{on }\ \Omega\times \{0\},\label{initial}
\end{eqnarray}
where $\displaystyle \Delta = \sum_{j=1}^{N+1} \frac {\partial^2}{\partial x_j^2}$,  or the Cauchy  problem:
\begin{equation}
\label{cauchy}
\partial_t U=\Delta \phi(U)\ \mbox{ in }\ \mathbb R^{N+1} \times (0, +\infty)\quad\mbox{ and }\ U = \chi_{\Omega^c}\ \mbox{ on }\ \mathbb R^{N+1}\times \{0\};
\end{equation}
here $\chi_{\Omega^c}$ denotes the characteristic function of the set $\Omega^c = \mathbb R^{N+1} \setminus \Omega$. Note that
the uniqueness of the solution of either problem \eqref{diffusion}-\eqref{initial} or problem \eqref{cauchy} follows from the 
comparison principle (see \cite[Theorem A.1, p. 253]{MSjde2}).
We consider the solution $U \in C^{2,1}(\Omega\times(0,+\infty)) \cap L^\infty(\Omega \times (0,+\infty))\cap C(\overline{\Omega} \times (0,+\infty)) $ such that $U(\cdot,t) \to 0$ in $L^1_{loc}(\Omega)$ as $t \to 0^+$ for problem \eqref{diffusion}-\eqref{initial}.  For problem  \eqref{cauchy},  we consider the solution $U \in C^{2,1}(\mathbb R^{N+1} \times (0,+\infty))\cap L^\infty(\mathbb R^{N+1} \times (0,+\infty))$ such that $U(\cdot,t) \to \chi_{\Omega^c}(\cdot)$ in $L^1_{loc}(\mathbb R^{N+1})$ as $t \to 0^+$.

By  the strong comparison principle, we know that
\begin{equation}
\label{strict bounds}
0 < U < 1\ \mbox{ and } \frac {\partial U}{\partial x_{N+1}} < 0\  \mbox{either  in } \Omega \times (0, +\infty) \mbox{ or in } \mathbb R^{N+1} \times(0,+\infty).
\end{equation}
The profile of $U$ as $t\to 0^+$ is controlled by the function $\Phi$ defined by
\begin{equation}
\label{definition of Phi}
\Phi(s) = \int_1^s \frac {\phi^\prime(\xi)}{\xi} d\xi\quad\mbox{ for } s > 0.
\end{equation}
In fact, in \cite[Theorem 2.1 and Remark 2.2, p. 239]{MSjde2} (see also \cite[Theorem 1.1 and Theorem 4.1, p. 940 and p. 947]{MSpoincare}) it is shown that, 
if $U$ is the solution of either problem {\rm (\ref{diffusion})-(\ref{initial})} or problem {\rm (\ref{cauchy})}, then
\begin{equation}
\label{varadhan formula}
-4t\Phi(U(X, t)) \to d(X)^2 \ \mbox{ as $t \to 0^+$\ uniformly on every compact subset of } \Omega.
\end{equation}
Here, $d = d(X)$ is the distance function:
\begin{equation}
\label{definition of distance}
d(X) = \mbox{ dist}(X, \partial\Omega)\quad\mbox{ for }  X=(x,x_{N+1}) \in \Omega.
\end{equation}
Formula \eqref{varadhan formula} is regarded as a nonlinear version of one obtained by Varadhan \cite{Va}.

A hypersurface $\Gamma$ in $\Omega$ is said to be a {\it stationary level surface} of $U$ ({\it stationary isothermic surface} of $U$ when $\phi(s) \equiv s$)
if at each time $t$ the solution $U$ remains constant on $\Gamma$ (a constant depending on $t$). Hence it follows from \eqref{varadhan formula} that there exists $R > 0$ such that
\begin{equation}
\label{constant distance}
d(X) = R\ \mbox{ for every } X \in \Gamma,
\end{equation}
provided $\Gamma$ is a stationary level surface of $U$.
The following theorem characterizes the boundary $\partial\Omega$ 
in such a way that $U$ has a stationary level surface $\Gamma$ in $\Omega.$

\begin{theorem}
\label{th:generalized BCN}
Let $U$ be the solution of either problem {\rm (\ref{diffusion})-(\ref{initial})} or problem {\rm (\ref{cauchy})}.
Assume that there exists a basis $\{y^1, y^2, \dots, y^N \} \subset \mathbb R^{N}$ such that for every $j = 1, \dots, N$ the function $f(x+y^{j}) - f(x)$ has either a maximum or a minimum in $\mathbb R^{N}$.  Suppose that $U$ has a stationary level surface $\Gamma$ in $\Omega$. Then $f$ is affine and $\partial\Omega$ must be a hyperplane.
\end{theorem}
\begin{remark}{\rm 
\label{concern Th 1.1} In order to prove Theorem \ref{th:generalized BCN}, we shall also use the sliding method due to Berestycki, Caffarelli, and Nirenberg \cite{BCN}. In \cite[Theorem 2.3 and Remark 2.4, p. 240]{MSjde2}, instead of the assumption on $f$, it is assumed that for each $y\in \mathbb R^{N}$ there exists $h(y) \in \mathbb R$ such that
\begin{equation}
\label{BCNmodified}
\lim_{|x| \to \infty} \left[ f(x + y) - f(x) \right] = h(y),
\end{equation}
which implies the assumption on $f$ in Theorem \ref{th:generalized BCN}. The condition \eqref{BCNmodified} is a modified version of \cite[(7.2), p. 1108]{BCN}, 
in which $h(y)$ is supposed identically zero. When $N=1$, $f(x) = ax+b+\sin x\ (a, b \in \mathbb R)$ satisfies the assumption on $f$ in Theorem \ref{th:generalized BCN}, but it does not satisfy \eqref{BCNmodified} provided $\frac y{2\pi}$ is not an integer. Another $f(x) = ax +b +\sin x\tan^{-1}x \ (a, b \in \mathbb R)$ does not satisfy the assumption, but it is Lipschitz continuous on $\mathbb R$. }
\end{remark}

Let us consider the case where $\phi(s) \equiv s$, that is, that of the heat equation, in particular. 
The following theorem characterizes the boundary $\partial\Omega$ 
in such a way that the caloric function $U$ has a stationary {\it isothermic} surface in $\Omega.$

\begin{theorem}
\label{th: caloric function}
Let $\phi(s) \equiv s$ and let $U$ be the solution of either problem {\rm (\ref{diffusion})-(\ref{initial})} or problem {\rm (\ref{cauchy})}. Assume that $U$ has a stationary isothermic surface $\Gamma$ in $\Omega$. Then $f$ is affine and $\partial\Omega$ must be a hyperplane, if either $N \le 2$ or  $\{ |f(x)-f(y)|: |x-y| \le 1 \}$ is bounded.
\end{theorem}

\begin{remark}{\rm 
\label{concern Th 1.3} When $f$ is Lipschitz continuous in $\mathbb R^N$ and $\Omega$ satisfies the uniform exterior sphere condition, this theorem was proved in \cite[Theorem 1.1 (ii), p. 1113]{MSjde1}. By combining \cite[Lemma 3.1]{MSmmas} with \cite[Theorem 1.1, p. 887]{S}, we see that the assumption that $\Omega$ satisfies the uniform exterior sphere condition is not needed. Also, the Lipschitz continuity of $f$ can be replaced by the uniform continuity of $f$, because of Professor Hitoshi Ishii's suggestion. Namely, by essentially the same proof as in \cite{S}, it can be shown that \cite[Theorem 1.1, p. 887]{S} holds even if the Lipschitz continuity is replaced by the uniform continuity. {\it Here, the advantage of {\rm Theorem \ref{th: caloric function}} is that we do not need to assume any uniform continuity of $f$.}
}
\end{remark}

Let $F= F(s)$ be a $C^1$ symmetric and concave function on the positive cone $\Lambda$ given by
$$
\Lambda =\{ s = (s_1,\cdots, s_N) \in \mathbb R^N  :  \min\limits_{1\le j\le N}s_j > 0 \},
$$
where $N \ge 1$. Assume that $F$  satisfies
\begin{equation}
\label{monotonicity}
 F_{s_j} \left( = \frac {\partial F}{\partial s_j} \right) > 0\ \mbox{ for all } j = 1, \cdots, N\ \mbox{ in }\Lambda.
\end{equation}
Define $G = G(s)$ by
\begin{equation}
\label{definition of G}
G(s) = F(1/{s_1}, \cdots, 1/{s_N})\ \mbox{ for } s \in \Lambda.
\end{equation}
Assume that $G$ is convex in $\Lambda$. Such a class of functions $F$ is dealt with in \cite{A, S}.
Related to Theorems \ref{th:generalized BCN} and \ref{th: caloric function}, for $f \in C(\mathbb R^N)$ we consider the domain $\Omega$ given by \eqref{upper-half domain}.
Consider the entire graph $\partial\Omega = \{ (x,f(x)) \in \mathbb R^{N+1} :  x \in \mathbb R^N \}$  in $\mathbb R^{N+1}$   of  $f$. Let $\kappa_1(x), \cdots, \kappa_N(x)$ be the principal curvatures of $\partial\Omega$ with respect to the upward unit normal vector to $\partial\Omega$ at $(x,f(x))$ for $x \in\mathbb R^N$. For each $R > 0$, we  introduce a function $g \in C(\mathbb R^N)$ defined by
\begin{equation}
\label{sup convolution Jensen}
g(x) = \sup_{|x-y| \le R}\left\{ f(y) +\sqrt{ R^2-|x-y|^2} \right\}\ \mbox{ for every } x \in \mathbb R^{N}.
\end{equation}
Then we have
\begin{equation}
\label{distance R property}
\{ (x, g(x)) \in \mathbb R^{N+1} : x \in \mathbb R^N \} = \{ X \in \mathbb R^{N+1} : d(X) = R \}\ (=\Gamma). 
\end{equation}
Moreover, let us introduce a function $f^* \in C(\mathbb R^N)$ defined by
\begin{equation}
\label{inf convolution Jensen}
f^*(x) = \inf_{|x-y| \le R}\left\{ g(y) -\sqrt{ R^2-|x-y|^2} \right\}\ \mbox{ for every } x \in \mathbb R^{N}.
\end{equation}
Then, by setting 
\begin{equation}
\label{upper half divided by Gamma}
D = \{ X = (x, x_{N+1}) \in \mathbb R^{N+1} : x_{N+1} > g(x)\ \}, 
\end{equation}
 we notice the following:
\begin{eqnarray}
&& \{ (x, f^*(x)) \in \mathbb R^{N+1} : x \in \mathbb R^N \} = \{ X \in \mathbb R^{N+1} : \mbox{ dist}(X, \overline{D}) = R \}, \label{converse distance property}
\\
&& f(x) \le f^*(x) \ \mbox{ for every }\ x \in \mathbb R^N. \label{a geometric inequality}
\end{eqnarray}
The third theorem gives a Liouville-type theorem for some Weingarten hypersurfaces in the viscosity sense.
\begin{theorem}
\label{th: viscosity}
Suppose that there exist two real constants $R > 0$ and $c$ such that  $f\in C(\mathbb R^N)$ satisfies 
  in the viscosity sense
 \begin{equation}
 \label{viscosity solutions}
F(1-R\kappa_1,\cdots,1-R\kappa_N) = c \  \mbox{ in }\  \mathbb R^N,
\end{equation}
and moreover suppose that the equality holds in \eqref{a geometric inequality}, that is, 
\begin{equation}
\label{the geometric equality}
f(x) = f^*(x)\ \mbox{ for every } \ x \in \mathbb R^N, 
\end{equation}
where $f^* = f^*(x)$ is defined by \eqref{inf convolution Jensen}.
Then, $c = F(1,\cdots,1)$ and $f$ is an affine function, that is, $\partial\Omega$ must be a hyperplane, provided  $\{ |f(x)-f(y)|: |x-y| \le 1 \}$ is bounded.
\end{theorem}
\begin{remark}{\rm 
\label{concern Th 1.5} The case where $F(s) = \left(\prod_{j=1}^N s_j\right)^{1/N}$ or  $F(s) = \sum\limits_{j=1}^N\log s_j$ is related to Theorem \ref{th: caloric function}. The assumption \eqref{the geometric equality}, that is, 
\begin{equation}
\label{inf convolution Jensen2}
f(x) = \inf_{|x-y|\le R}\left\{ g(y) -\sqrt{R^2- |x-y|^{2}}\right\}\ \mbox{ for every } x \in \mathbb R^N,
\end{equation}
implies that 
\begin{equation}
\label{principal curvature bound}
\max_{1\le j \le N} \kappa_j \le \frac 1R\ \mbox{ in }\ \mathbb R^N
\end{equation}
holds in the viscosity sense, because \eqref{the geometric equality} yields that for every point $X \in \partial\Omega$  there exists an open ball $B_R(Y)$ with radius $R$ and centered at $Y \in \Gamma$ satisfying
\begin{equation}
\label{useful ball property1}
X \in \partial B_R(Y)\ \mbox{ and } B_R(Y) \subset \Omega.
\end{equation}
\eqref{principal curvature bound} is one of main assumptions of \cite[Theorem 1.1, p. 887]{S}. Namely, boundedness of $\{ |f(x)-f(y)|: |x-y| \le 1 \}$  is much weaker than Lipschitz continuity of $f$, but \eqref{the geometric equality} is stronger than \eqref{principal curvature bound}.  Also, \eqref{the geometric equality} is satisfied by every classical $C^2$ solution $f$ of \eqref{viscosity solutions} having the strict inequality in \eqref{principal curvature bound}, because of the implicit function theorem. } 
\end{remark}

The present paper is organized as follows. In Section \ref{Proof of 1st theorem},  we prove Theorem \ref{th:generalized BCN} by using the sliding method due to Berestycki, Caffarelli, and Nirenberg \cite{BCN}. 
In Section \ref{Proof of 2nd theorem},  we prove Theorem \ref{th: caloric function} with the aid of the theory of viscosity solutions. We follow the proof of  \cite[Theorem 1.1, p. 887]{S}  basically, but we here need a key lemma (see Lemma \ref{le: gradient estimates}) which gives new gradient estimates for $f$ and $g$, because we do not assume any uniform continuity of $f$.
Section \ref{Proof of 3rd theorem} is devoted to a proof of Theorem \ref{th: viscosity}, where gradient estimates for $f$ and $g$ are replaced by  Lipschitz constant estimates for $f$ and $g$ (see Lemma \ref{le: Lipschitz estimates}). In Section \ref{concluding remarks}, we give a Bernstein-type theorem for some $C^{2}$ Weingarten hypersurfaces (see Theorem \ref{th: classical Bernstein theorem for Weingarten}) as a remark on Theorem \ref{th: viscosity}.


\setcounter{equation}{0}
\setcounter{theorem}{0}

\section{Proof of Theorem \ref{th:generalized BCN}}
\label{Proof of 1st theorem}

Since $\Gamma$ is a stationary level surface of $U$, 
it follows from \eqref{varadhan formula}, \eqref{strict bounds} and the implicit function theorem that 
there exist a number $R > 0$ and a function $g \in C^2(\mathbb R^N)$ such that both \eqref{sup convolution Jensen} and \eqref{distance R property} hold.

Conversely, let $\nu(y)$ denote the upward unit normal vector to $\Gamma$ at $(y, g(y)) \in \Gamma.$  The facts that $g$ is smooth, $\partial\Omega$ is a graph, and $(y, g(y)) - R\nu(y) \in \partial\Omega$ for every $y\in\mathbb R^{N}$, imply that \eqref{the geometric equality}, \eqref{inf convolution Jensen}, and \eqref{converse distance property} hold, namely, both \eqref{inf convolution Jensen2} and \eqref{converse distance property} where $f^*$ is replaced by $f$ hold. Hence, we have in particular
\begin{equation}
\label{rep of partial Omega}
\partial\Omega = \{ (x, f(x)) \in \mathbb R^{N+1} : x \in \mathbb R^N \} = \{ X \in \mathbb R^{N+1} : \mbox{ dist}(X, \overline{D}) = R\},
\end{equation}
where $D$ is given by \eqref{upper half divided by Gamma}.
Thus, it follows from \eqref{rep of partial Omega}  that for every $X \in \partial\Omega$ there exists $Y \in \Gamma$ satisfying 
\begin{equation}
\label{useful ball touching partial Omega}
X \in \partial B_R(Y)\ \mbox{ and }\ B_R(Y) \subset \Omega.
\end{equation}
Choose $j$ arbitrarily. By the assumption of Theorem \ref{th:generalized BCN}, the function $f(x+y^{j}) - f(x)$ has either a maximum or a minimum in $\mathbb R^{N}$. Since the proof below is similar, say $f(x+y^{j}) - f(x)$ has  a maximum $M$ in $\mathbb R^{N}$. Then there exists $x_0 \in \mathbb R^N$ such that
\begin{equation}
\label{maximum is attained}
f(x +y^j)- f(x) \le M = f(x_0+y^j) -f(x_0)\ \mbox{ for every } x \in \mathbb R^N.
\end{equation}
Let us use the sliding method due to Berestycki, Caffarelli, and Nirenberg \cite{BCN}. We set
$$
\Omega_{y^j, M} = \{ (x, x_{N+1}) \in \mathbb R^{N+1} : (x+y^j, x_{N+1} + M) \in \Omega \}.
$$
Then we have
\begin{eqnarray*}
&& f(x+y^j) -M \le f(x)\ \mbox{ for every } x \in \mathbb R^N,
\\
&& \Omega_{y^j, M} \supset \Omega\ \mbox{ and } (x_0, f(x_0)) \in \partial\Omega \cap \partial \Omega_{y^j, M}.
\end{eqnarray*}
Suppose that $\Omega_{y^j, M} \supsetneqq \Omega$. Then, by the strong comparison principle we have
\begin{equation}
\label{strict comparison with translate}
 U(x+y^j, x_{N+1}+M, t) < U(X,t) \mbox{ for every } (X, t) = (x, x_{N+1}, t) \in \Omega \times (0,  +\infty).
\end{equation}
On the other hand, since $(x_0, f(x_0)) \in \partial\Omega \cap \partial \Omega_{y^j, M}$ and $\Omega_{y^j, M} \supset \Omega$, it follows from \eqref{useful ball touching partial Omega}  that there exists $Y_{0}=(y_{0}, g(y_{0})) \in \Gamma$ satisfying
$$
(x_0, f(x_0)) \in \partial B_R(Y_{0})\ \mbox{ and }\ B_R(Y_{0}) \subset \Omega \subset \Omega_{y^j, M}.
$$
Hence, since $\Gamma = \{ X \in \mathbb R^{N+1} : d(X) = R \}$ and $\Gamma$ is a stationary level surface of $U$, we have
$$
U(y_{0}+y^j, g(y_{0})+M, t) = U(Y_{0},t) \ \mbox{ for every } t > 0,
$$
which contradicts \eqref{strict comparison with translate}. Thus, we get $\Omega_{y^j, M} = \Omega$, that is,
$$
f(x +y^j) - M = f(x)\ \mbox{ for every } x \in \mathbb R^N.
$$
Therefore we conclude that there exist $a_1, \dots, a_N \in \mathbb R$ satisfying
\begin{equation}
\label{first stage functional equations}
f(x +y^j) = f(x) + a_j\ \mbox{ for every } x \in \mathbb R^N \mbox{ and for } j=1, \dots, N,
\end{equation}
since $j$ is chosen arbitrarily. Since $f$ is continuous on $\mathbb R^N$ and $\{y^1, y^2, \dots, y^N \}$ is a basis of $ \mathbb R^N$, we can solve \eqref{first stage functional equations} as a system of functional equations and conclude that $f(x)$ is determined by its values on $E = \{ \sum_{j=1}^N\beta_j y^j\in \mathbb R^N : 0 \le \beta_j < 1,\ j= 1, \dots, N \}$. Indeed, if $x = \sum_{j=1}^N(r_j+\beta_j)y^j$ for $r =(r_1, \dots, r_N) \in \mathbb Z^N$ and  $\beta =(\beta_1, \dots, \beta_N) \in [0, 1)^N$, then $f(x) = f( \sum_{j=1}^N\beta_j y^j) + \sum_{j=1}^N r_j a_j$. Moreover, this property of $f$ implies that for every $y \in \mathbb R^N$ the function
$f(x+y) - f(x)$ has either a maximum or a minimum on $\mathbb R^N$. Thus, by employing the sliding method again, we get
\begin{equation}
\label{second stage functional equations}
f(x +y) - f(x) = f(z+y) - f(z)\ \mbox{ for every } x, y, z \in \mathbb R^N.
\end{equation}
Since $f$ is continuous  on $\mathbb R^N$, we solve \eqref{second stage functional equations} as a system of functional equations and conclude that $f$ is affine.
This completes the proof of Theorem \ref{th:generalized BCN}.


\setcounter{equation}{0}
\setcounter{theorem}{0}

\section{Proof of Theorem \ref{th: caloric function}}
\label{Proof of 2nd theorem}
Note that $U$ is real analytic in $x$, since $U$ satisfies the heat equation.
Since $\Gamma$ is a stationary isothermic surface of $U$, it follows from \eqref{strict bounds} and the implicit function theorem that $\Gamma$ is the graph of a real analytic function $g= g(x)$ for $x \in \mathbb R^N$.
Let us first quote an important lemma from \cite[Lemma 3.1]{MSmmas}.
We can use this lemma, since $\partial\Omega =\partial\left(\mathbb R^{N+1}\setminus\overline{\Omega}\right)$, $\Gamma$ is already real analytic and $\Gamma = \partial D$ where $D$ is given by \eqref{upper half divided by Gamma}. The interior cone condition of $D \mbox{ in the lemma })$ with respect to $\Gamma$ is of course satisfied, but in \cite{MSmmas} it is used only to show that $\Gamma$ is smooth.

\begin{lemma}[\cite{MSmmas}]
\label{le:preliminary} The following assertions hold:
\begin{itemize}
\item[\rm (1)]\ There exists a number $R > 0$ such that $d(X) = R$ for every $X \in \Gamma$; 
\item[\rm (2)]\ $\Gamma$ is a real analytic hypersurface;
\item[\rm (3)]\ $\partial\Omega$ is also a real analytic hypersurface,
such that the mapping $\partial\Omega \ni (x, f(x)) \mapsto Y(x, f(x)) \equiv (x, f(x)) + R\nu(x) \in \Gamma,$ 
where $\nu(x)$ is the upward unit normal vector to $\partial\Omega$ at $(x,f(x)) \in \partial\Omega$, is a diffeomorphism; 
in particular, $\partial\Omega$ and $\Gamma$ are parallel hypersurfaces at distance $R$;
\item[\rm (4)]\ it holds that
\begin{equation}
\label{strict bounds of curvatures}
 \max_{1\le j \le N}\kappa_j(x) < \frac 1R\ \mbox{ for every } x \in \mathbb R^{N},
\end{equation}
where $\kappa_1(x), \dots, \kappa_{N}(x)$ are the principal curvatures of $\partial\Omega$ at $(x,f(x)) \in \partial\Omega$ with respect to the upward unit normal vector to $\partial\Omega$;
\item[\rm (5)]\ there exists a number $c > 0$ such that
\begin{equation}
\label{monge-ampere}
\prod_{j=1}^{N} \left(1-R\kappa_j(x)\right) = c\quad\mbox{ for every } x \in \mathbb R^{N}.
\end{equation}
\end{itemize}
\end{lemma}
Note that (1) follows from \eqref{varadhan formula} and (2) follows simply from the implicit function theorem.
When $N = 1$, (5) implies the conclusion of Theorem \ref{th: caloric function}, since the curvature of the curve $\partial\Omega$ is constant. Let $N \ge 2$. With the aid of Lemma \ref{le:preliminary}, applying \cite[Lemmas 4.2 and 4.3, p. 891 and p. 892]{S} to $F(s) = \left(\prod_{j=1}^N s_j\right)^{1/N}$ yields the following lemma.

\begin{lemma}
\label{le:mean curvature inequality}
$c = 1$ and $H_{\partial\Omega}\le 0\le H_{\Gamma}$ in $\mathbb R^{N}$, where $H_{\partial\Omega}$ {\rm (resp. $H_{\Gamma}$)} denotes the mean curvature of $\partial\Omega$ {\rm (resp. $\Gamma$)} with respect to the upward unit normal vector to $\partial\Omega$ {\rm (resp. $\Gamma$)}.
\end{lemma} 

When $N =2$, by setting 
\begin{equation}
\label{middle parallel surface}
\Gamma^{*} = \left\{ X \in \Omega : d(X) = \frac R2 \right\},
\end{equation}
the fact that $c = 1$ implies that $\Gamma^{*}$ is an entire minimal graph over $\mathbb R^{2}$. Therefore, by the Bernstein's theorem for the minimal surface equation, $\Gamma^{*}$ must be a hyperplane as in \cite{MSjde1}.
(See \cite{GT, Gi} for the Bernstein's theorem, and for more general setting see also Theorem \ref{th: classical Bernstein theorem for Weingarten}  in Section \ref{concluding remarks} in the present paper.)
Thus it remains to consider the case where $N\ge 3$ and $\{ |f(x)-f(y)|: |x-y| \le 1 \}$ is bounded.

On the other hand, (3) of Lemma \ref{le:preliminary} gives us the following geometric property.
\begin{lemma}
\label{le:useful geometry}
The following two assertions hold:
\begin{itemize}
\item[\rm (i)]\ For every $Y \in \Gamma$ there exists $X \in \partial\Omega$ such that
$
Y \in \partial B_R(X) \mbox{ and } B_R(X) \subset \mathbb R^{N+1}\setminus \overline{D}.
$
\item[\rm (ii)]\ For every $X \in \partial\Omega$ there exists $Y \in \Gamma$ such that
$
X \in \partial B_{R}(Y) \mbox{ and } B_{R}(Y) \subset \Omega.
$
\end{itemize}
\end{lemma}
Recall that $f$ and $g$ have the relationship, \eqref{sup convolution Jensen} and \eqref{inf convolution Jensen2}.
Since $\{ |f(x)-f(y)|: |x-y| \le 1 \}$ is bounded, we see that $\{ |g(x)-g(y)|: |x-y| \le 1 \}$ is also bounded. 
By Lemma \ref{le:mean curvature inequality} we have
\begin{equation}
\label{mean curvature inequality}
\mathcal M(f) \le 0 \le \mathcal M(g) \equiv \mbox{ div}\left(\frac {\nabla g}{\sqrt{1+|\nabla g|^{2}}}\right)\ \mbox{ in } \mathbb R^{N}.
\end{equation}
Let $B_{n}= \{ x \in \mathbb R^{N} : |x| < n \}$ for  $n \in \mathbb N$. Then, by \cite[Theorem 16.9, pp. 407--408]{GT}, for each $n \in \mathbb N$, there exist two functions $f_{n}, g_{n} \in C^{2}(B_{n}) \cap C(\overline{B_{n}})$ solving
\begin{eqnarray*}
&&\mathcal M(f_{n}) = \mathcal M(g_{n}) = 0\ \mbox{ in } B_{n},
\\
&& f_{n}= f\ \mbox{ and } g_{n}=g\ \mbox{ on } \partial B_{n}.
\end{eqnarray*}
Hence it follows from the comparison principle that for each $n\in\mathbb N$ there exists $z_{n}\in \partial B_{n}$ such that
\begin{equation}
\label{approximating sequence}
 f_{n+1} \le f_{n}\le f < g \le g_{n}\le g_{n+1}\ \mbox{ and } g_{n}-f_{n} \le g(z_{n})- f(z_{n})\ \mbox{ in } B_{n}.
\end{equation}
Since $\{ |f(x)-f(y)|: |x-y| \le 1 \}$ is bounded, it follows from  \eqref{sup convolution Jensen}   that $g-f$ is bounded in $\mathbb R^N$ and hence with the aid of  \eqref{approximating sequence} there exists a constant $C_* > 0$ satisfying
\begin{equation}
\label{boundedness of  f and g}
g-C_* \le f_n \le f \ \mbox{ and }\ g \le g_n \le f + C_*\ \mbox{ in } B_n \ \mbox{ for every } n \in \mathbb N.
\end{equation}
Thus, since both $\{ |f(x)-f(y)|: |x-y| \le 1 \}$ and $\{ |g(x)-g(y)|: |x-y| \le 1 \}$ are bounded, by using the interior estimates for the minimal surface equation (see \cite[Corollary 16.7, p. 407]{GT}) with the aid of   \eqref{boundedness of  f and g} and the monotonicity with $n$ in \eqref{approximating sequence}, we proceed as in \cite[pp. 893--894]{S} to see that there exist two functions $f_{\infty},  g_{\infty} \in C^{\infty}(\mathbb R^{N})$ satisfying
\begin{eqnarray*}
&& \mathcal M(f_{\infty}) = \mathcal M(g_{\infty}) = 0\ \mbox{ in } \mathbb R^{N},
\\
&& |\nabla f_{\infty}| \mbox{ and } |\nabla g_{\infty}| \mbox{ are bounded on }\mathbb R^{N},
\\
&& f_{n}\to f_{\infty}\ \mbox{ and } g_{n}\to g_{\infty}\ \mbox{ as $n \to \infty$ uniformly on every compact set in } \mathbb R^{N}.
\end{eqnarray*}
Then it follows from Moser's theorem \cite[Corollary, p. 591]{Mo} that both $f_{\infty}$ and $g_{\infty}$ are affine and hence the graph of $f_{\infty}$ is parallel to that of $g_{\infty}$ because $f_{\infty}\le g_{\infty}$ in $\mathbb R^{N}$. Thus there exists $\eta \in \mathbb R^{N}$ satisfying
\begin{equation}
\label{two affine functions}
f_{\infty}(x) = \eta\cdot x + f_{\infty}(0)\ \mbox{ and }\ g_{\infty}(x) = \eta\cdot x + g_{\infty}(0)\ \mbox{ for every } x \in \mathbb R^{N}.
\end{equation}
Moreover we have
\begin{eqnarray}
&& f_{\infty}\le f < g \le g_{\infty}\ \mbox{ in } \mathbb R^{N}, \label{order preserved}
\\
&& f(z_{n})- f_{\infty}(z_{n})\ \mbox{ and }\ g_{\infty}(z_{n})-g(z_{n})\ \to 0\ \mbox{ as } n \to \infty.\label{zero sequences}
\end{eqnarray}
Indeed, \eqref{order preserved} follows from \eqref{approximating sequence}. Observe that for each $n \in\mathbb N$
\begin{eqnarray*}
g_{n}(0) -f_{n}(0) &\le& g(z_{n}) - f(z_{n}) \le g_{n+1}(z_{n}) - f_{n+1}(z_{n}) 
\\
&\le& g(z_{n+1}) - f(z_{n+1}) \le
g_{\infty}(z_{n+1})-f_{\infty}(z_{n+1}) = g_{\infty}(0) - f_{\infty}(0).
\end{eqnarray*}
Hence letting $n \to \infty$ yields that
$
g(z_{n}) - f(z_{n}) \to g_{\infty}(0) - f_{\infty}(0)\ \mbox{ as } n \to \infty.
$
Thus as $n \to \infty$
$$
(f(z_{n})- f_{\infty}(z_{n})) + (g_{\infty}(z_{n})-g(z_{n})) = (g_{\infty}(0) - f_{\infty}(0)) - (g(z_{n}) -f(z_{n})) \to 0,
$$
which gives \eqref{zero sequences}.

It suffices to show that $f \equiv f_{\infty}$ and $g \equiv g_{\infty}$. Lemma \ref{le:useful geometry} yields the following key lemma.
\begin{lemma}[gradient estimates]
\label{le: gradient estimates}
There exist three constants $\varepsilon_{0} > 0, \delta_{0} > 0,$ and $C_{0}> 0$ such that
\begin{itemize}
\item[\rm (1)]\ if $z \in \mathbb R^{N}$ and $(0 \le) g_{\infty}(z) - g(z) \le \varepsilon_{0}$, then
$\displaystyle
\sup_{|y-z|\le\delta_{0}}|\nabla g(y)| \le C_{0}.
$
\item[\rm (2)]\ if $z \in \mathbb R^{N}$ and $(0 \le) f(z) - f_{\infty}(z) \le \varepsilon_{0}$, then
$\displaystyle
\sup_{|x-z|\le\delta_{0}}|\nabla f(x)| \le C_{0}.
$
\end{itemize}
\end{lemma}

\noindent
{\it Proof. } (i) of Lemma \ref{le:useful geometry} yields (1) and (ii) of Lemma \ref{le:useful geometry} yields (2), respectively. Let us show (1). Recall that $g_{\infty}$ is affine and $\nabla g_{\infty}\equiv\eta$. Denote by  $\mathcal H$ the hyperplane given by the graph of $g_{\infty}$. Then $\displaystyle \frac {(-\eta, 1)}{\sqrt{1 + |\eta|^{2}}}$ is the upward unit normal vector to $\mathcal H$. By (i) of Lemma \ref{le:useful geometry}, for every $Y=(y,g(y)) \in \Gamma$ there exists $X=(x,f(x)) \in \partial\Omega$ such that the ball $B_{R}(X)$ touching $\Gamma$ from below at $Y  \in \Gamma$ must be below  $\mathcal H$.  Hence, \begin{equation}
\label{normal approximation}
\mbox{if $Y$ is sufficiently close to $\mathcal H$, then }
\frac {Y-X}R \mbox{ is sufficiently close to } \frac {(-\eta, 1)}{\sqrt{1 + |\eta|^{2}}}. 
\end{equation}
Namely, for every $\mu > 0$ there exists $\lambda >0$ such that, if $(0 \le) g_{\infty}(y) - g(y) \le \lambda$, then 
\begin{equation}
\label{close in quantity}
\left|\frac{y-x}{R}- \frac {-\eta}{\sqrt{1 + |\eta|^{2}}}\right|^{2}+\left(\frac{g(y)-f(x)}{R}- \frac {1}{\sqrt{1 + |\eta|^{2}}}\right)^{2}< \mu^{2}.
\end{equation}
Of course, at the touching point $Y$, $\nabla g(y)$ equals the gradient of $f(x) + \sqrt{R^{2}-|y-x|^{2}}$ with respect to $y$, that is, 
\begin{equation}
\label{equal gradients}
\nabla g(y) = -\frac {y-x}{\sqrt{R^2-|y-x|^{2}}}.
\end{equation}
 On the other hand, if a point $(z,g(z)) \in \Gamma$ is sufficiently close to $\mathcal H$, then by \eqref{normal approximation} there exists a uniform neighborhood $\mathcal N_{z}$ of $z$ in $\mathbb R^{N}$ such that
every point $Y =(y,g(y)) \in \Gamma$ with $y \in \mathcal N_{z}$ is sufficiently close to $\mathcal H$. Namely, for every $\lambda > 0$ there exist $\varepsilon>0$ and $\delta > 0$ such that, if $(0 \le) g_{\infty}(z) - g(z) \le \varepsilon$ and $|y-z| < \delta$, then $(0 \le) g_{\infty}(y) - g(y) \le \lambda$. Thus, combining this fact with \eqref{close in quantity} and \eqref{equal gradients} yields (1). (2) is similar. \qed

The last lemma is
\begin{lemma}
\label{le:convergence}
The following two assertions hold:
\begin{itemize}
\item[\rm (i)]\ $g_{\infty}(x+z_{n}) - g(x+z_{n}) \to 0$ as $n \to \infty$ uniformly on every compact set in $\mathbb R^{N}$.
\item[\rm (ii)]\ $f(x+z_{n}) - f_{\infty}(x+z_{n}) \to 0$ as $n \to \infty$ uniformly on every compact set in $\mathbb R^{N}$.
\end{itemize}
\end{lemma}
This lemma implies the conclusion of Theorem \ref{th: caloric function}. Indeed, 
in view of \eqref{order preserved} and Lemma \ref{le:useful geometry}, Lemma \ref{le:convergence} yields that 
the graphs of $g_{\infty}$ and $f_{\infty}$ are parallel hyperplanes at distance $R$. This means that 
$f \equiv f_{\infty}$ and $g \equiv g_{\infty}$. Thus it remains to prove Lemma \ref{le:convergence}.

\vskip 2ex
\noindent
{\it Proof of Lemma \ref{le:convergence}. } Since (ii) is similar to (i), let us show (i). Set 
$$
G_{n}(x) = g(x + z_{n})-g(z_{n})\ \mbox{ for } x \in \mathbb R^{N} \mbox{ and } n \in \mathbb N.
$$
Then $G_{n}(0) = 0$ for every $n \in \mathbb N$. Since by \eqref{zero sequences} $g_{\infty}(z_{n})-g(z_{n}) \to 0$ as $n \to \infty$, it follows from (1) of Lemma \ref{le: gradient estimates} that there exists $N_{0}\in \mathbb N$ such that $\{ G_{n} : n \ge N_{0} \}$ is equicontinuous and bounded on $\overline{B_{\delta_{0}}(0)}\ (\subset \mathbb R^{N})$.
Arzela-Ascoli theorem gives us that there exist a subsequence $\{G_{n^{\prime}}\}$ and a function $G_{\infty} \in C(\overline{B_{\delta_{0}}(0)})$ such that
\begin{equation}
\label{uniform convergence of subsequence}
G_{n^{\prime}} \to G_{\infty}\ \mbox{ as } n \to \infty\ \mbox{ uniformly on } \overline{B_{\delta_{0}}(0)}.
\end{equation}
Notice that $G_{\infty}(0) = 0$. Since $\mathcal M(G_{n}) \ge 0$ in $\mathbb R^{N}$ by \eqref{mean curvature inequality}, we have that $\mathcal M(G_{\infty}) \ge 0$ in $B_{\delta_{0}}(0)$ in the viscosity sense. Observe that
\begin{eqnarray*}
&&G_{n^{\prime}}(x) =g(x + z_{n^{\prime}})-g(z_{n^{\prime}}) \le g_{\infty}(x + z_{n^{\prime}})-g(z_{n^{\prime}})
\\
&&\qquad= \left\{g_{\infty}(x + z_{n^{\prime}})- g_{\infty}(z_{n^{\prime}}) \right\} + \left\{g_{\infty}(z_{n^{\prime}}) -g(z_{n^{\prime}})\right\} = \eta\cdot x +  \left\{g_{\infty}(z_{n^{\prime}}) -g(z_{n^{\prime}})\right\}.
\end{eqnarray*}
Hence, by \eqref{zero sequences} and \eqref{uniform convergence of subsequence}, letting $n^{\prime}\to \infty$ yields
\begin{equation}
\label{under the hyperplane}
G_{\infty}(x) \le \eta\cdot x\ \mbox{ in } \overline{B_{\delta_{0}}(0)}.
\end{equation}
Therefore, since $\mathcal M( \eta\cdot x ) = 0 \le \mathcal M(G_{\infty})$ in $B_{\delta_{0}}(0)$ in the viscosity sense and $\eta\cdot 0 = 0 = G_{\infty}(0)$, by the strong comparison principle of Giga and Ohnuma \cite[Theorem 3.1, p. 173]{GO} we see that
$$
G_{\infty}(x) \equiv \eta\cdot x \ \mbox{ in } \overline{B_{\delta_{0}}(0)}.
$$
Thus $G_{\infty}$ is uniquely determined independently of the choice of the subsequence and therefore from \eqref{uniform convergence of subsequence} we conclude that 
\begin{equation}
\label{uniform convergence of original sequence}
G_{n}(x) \to \eta\cdot x \ \mbox{ as } n \to \infty \mbox{ uniformly on } \overline{B_{\delta_{0}}(0)}.
\end{equation}
Then, since 
\begin{eqnarray*}
g_{\infty}(x+z_{n}) - g(x+z_{n}) &=& \left\{ g_{\infty}(x+z_{n}) - g_{\infty}(z_{n}) \right\} - G_{n}(x) + \left\{ g_{\infty}(z_{n})- g(z_{n})\right\}
\\
&=& \eta\cdot x- G_{n}(x) + \left\{ g_{\infty}(z_{n})- g(z_{n})\right\},
\end{eqnarray*}
we get from \eqref{zero sequences} and \eqref{uniform convergence of original sequence}
\begin{equation}
\label{uniform convergence on the first ball}
g_{\infty}(x+z_{n}) - g(x+z_{n}) \to 0\ \mbox{ as } n \to \infty \mbox{ uniformly on } 
\overline{B_{\delta_{0}}(0)}.
\end{equation}
Moreover, by using (1) of Lemma \ref{le: gradient estimates} again for any point $z \in \partial B_{\delta_{0}}(0)$ and repeating the same argument as above, we see that \eqref{uniform convergence on the first ball} holds even if $B_{\delta_{0}}(0)$ is replaced by $B_{\frac 32\delta_{0}}(0)$. Thus, repeating this argument as many times as one wants yields conclusion (i). \qed

\begin{remark}{\rm 
\label{concern the case N=1}
For the proof of Theorem \ref{th: viscosity}, we give a remark for the case where $N=1$. Even when $N=1$, all the lemmas \ref{le:mean curvature inequality} - \ref{le:convergence} hold true. Indeed, 
when $N=1$,  $\displaystyle\mathcal M(g) =g^{\prime\prime}( 1+(g^\prime)^2)^{-\frac 32}$ in \eqref{mean curvature inequality}. Hence the graphs of $f_n$ and $g_n$ are line segments and without using Moser's theorem we can get two affine functions $f_\infty$ and $g_\infty$ in \eqref{two affine functions}.}
\end{remark}


\setcounter{equation}{0}
\setcounter{theorem}{0}

\section{Proof of Theorem \ref{th: viscosity}}
\label{Proof of 3rd theorem}

We follow the proof of Theorem \ref{th: caloric function}. By \cite[Lemmas 4.2 and 4.3, p. 891 and p. 892]{S}, we have instead of Lemma \ref{le:mean curvature inequality} 
\begin{lemma}
\label{le:mean curvature inequality viscosity sense}
$c = F(1, \cdots, 1)$ and $H_{\partial\Omega}\le 0\le H_{\Gamma}$ in $\mathbb R^{N}$ in the viscosity sense, where $H_{\partial\Omega}$ {\rm (resp. $H_{\Gamma}$)} denotes the mean curvature of $\partial\Omega$ {\rm (resp. $\Gamma$)} with respect to the upward unit normal vector to $\partial\Omega$ {\rm (resp. $\Gamma$)}.
\end{lemma} 

Also, in view of \eqref{sup convolution Jensen} and \eqref{inf convolution Jensen2} coming from \eqref{the geometric equality}, we see that Lemma \ref{le:useful geometry} also holds. Then proceeding as in the proof of Theorem \ref{th: caloric function} yields two affine functions $f_{\infty}$ and $g_{\infty}$ satisfying 
\eqref{two affine functions}, \eqref{order preserved}, and \eqref{zero sequences}. Hence, it suffices to show that $f \equiv f_{\infty}$ and $g \equiv g_{\infty}$. Lemma \ref{le:useful geometry} yields the following key lemma instead of Lemma \ref{le: gradient estimates}.

\begin{lemma}[Lipschitz constant estimates]
\label{le: Lipschitz estimates}
There exist three constants $\varepsilon_{0} > 0, \delta_{0} > 0,$ and $C_{0}> 0$ such that
\begin{itemize}
\item[\rm (1)]\ if $z \in \mathbb R^{N}$ and $(0 \le) g_{\infty}(z) - g(z) \le \varepsilon_{0}$, then
$\displaystyle
\sup_{x, y \in B_{\delta_{0}}(z), x \not=y} \frac {|g(x)-g(y)|}{|x-y|} \le C_{0}.
$
\item[\rm (2)]\ if $z \in \mathbb R^{N}$ and $(0 \le) f(z) - f_{\infty}(z) \le \varepsilon_{0}$, then
$\displaystyle
\sup_{x, y \in B_{\delta_{0}}(z), x \not=y} \frac {|f(x)-f(y)|}{|x-y|} \le C_{0}.
$
\end{itemize}
\end{lemma}

\noindent
{\it Proof. } We adjust the proof of Lemma \ref{le: gradient estimates} to this situation. (i) of Lemma \ref{le:useful geometry} yields (1) and (ii) of Lemma \ref{le:useful geometry} yields (2), respectively. Let us show (1). Recall that $g_{\infty}$ is affine and $\nabla g_{\infty}\equiv\eta$. Denote by  $\mathcal H$ the hyperplane given by the graph of $g_{\infty}$. Then $\displaystyle \frac {(-\eta, 1)}{\sqrt{1 + |\eta|^{2}}}$ is the upward unit normal vector to $\mathcal H$. By (i) of Lemma \ref{le:useful geometry},  for every $Y =(y,g(y)) \in \Gamma$ there exists $X=(x,f(x)) \in \partial\Omega$ such that the ball $B_{R}(X)$ touching $\Gamma$ from below at $Y  \in \Gamma$ must be below  $\mathcal H$.  Hence, 
\begin{equation}
\label{normal approximation viscosity}
\mbox{if $Y$ is sufficiently close to $\mathcal H$, then }
\frac {Y-X}R \mbox{ is sufficiently close to } \frac {(-\eta, 1)}{\sqrt{1 + |\eta|^{2}}}. 
\end{equation}
Namely, for every $\mu > 0$ there exists $\lambda >0$ such that, if $(0 \le) g_{\infty}(y) - g(y) \le \lambda$, then 
\begin{equation}
\label{close in quantity2}
\left|\frac{y-x}{R}- \frac {-\eta}{\sqrt{1 + |\eta|^{2}}}\right|^{2}+\left(\frac{g(y)-f(x)}{R}- \frac {1}{\sqrt{1 + |\eta|^{2}}}\right)^{2}< \mu^{2}.
\end{equation}

 On the other hand, if a point $(z,g(z)) \in \Gamma$ is sufficiently close to $\mathcal H$, then by \eqref{normal approximation viscosity} there exists a uniform neighborhood $\mathcal N_{z}$ of $z$ in $\mathbb R^{N}$ such that
every point $Y =(y,g(y)) \in \Gamma$ with $y \in \mathcal N_{z}$ is sufficiently close to $\mathcal H$. Namely, for every $\lambda > 0$ there exist $\varepsilon>0$ and $\delta > 0$ such that, if $(0 \le) g_{\infty}(z) - g(z) \le \varepsilon$ and $|y-z| < \delta$, then $(0 \le) g_{\infty}(y) - g(y) \le \lambda$.

Moreover, in view of \eqref{normal approximation viscosity}, by choosing $\frac \pi 2 > \theta > 0$ sufficiently small and introducing a cone $\mathcal V$ defined by
$$
\mathcal V = \{ \Xi=(\xi,\xi_{N+1})\in \mathbb R^{N+1} : \xi_{N+1} > |\Xi|\cos\theta\},
$$
we see that, if $Y \in \Gamma$ is sufficiently close to $\mathcal H$, then $\mathcal V + Y =\{ \Xi + Y\ : \Xi \in \mathcal V\} \subset D$, where $\mathcal V + Y$ is a cone with vertex $Y$. Here $D$ is given by \eqref{upper half divided by Gamma}. Indeed, if $\mathcal V + Y \not\subset D$, then there exists another point $\tilde{Y} (\not= Y) \in \Gamma\cap (\mathcal V+Y)$. However, in view of \eqref{normal approximation viscosity}, a ball $B_{R}(\tilde{X})$ touching $\Gamma$ from below at $\tilde{Y}$ might contain $Y$ since $\theta >0$ is small. This is a contradiction. Namely, if $Y =(y,g(y)) \in \Gamma$ with $y \in \mathcal N_{z}$, then, with the aid of (i) of Lemma \ref{le:useful geometry}, we must have 
$$
\mathcal V + Y  \subset D,\ B_R(X) \subset \mathbb R^{N+1}\setminus \overline{D},\ \mbox{ and } Y \in
\partial(\mathcal V+Y)\cap \partial B_R(X).
$$
This gives (1). (2) is similar. \qed

Hence, by using Lemma \ref{le: Lipschitz estimates} instead of Lemma \ref{le: gradient estimates}, we can proceed as in the proof of Theorem \ref{th: caloric function} to see that Lemma \ref{le:convergence} also holds.
Therefore, \eqref{order preserved}, Lemma \ref{le:useful geometry} and Lemma \ref{le:convergence} yield the conclusion of Theorem \ref{th: viscosity}.


\setcounter{equation}{0}
\setcounter{theorem}{0}

\section{Concluding remarks}
\label{concluding remarks}

When $N =2$, we have a Bernstein-type theorem for some $C^{2}$ Weingarten hypersurfaces related to Theorem \ref{th: viscosity}.
\begin{theorem}
\label{th: classical Bernstein theorem for Weingarten}
Suppose that there exist two real constants $R > 0$ and $c$ such that  $f\in C^{2}(\mathbb R^2)$ satisfies 
 \begin{equation}
 \label{classical solutions of Weingarten}
F(1-R\kappa_1, 1-R\kappa_2) = c \ \mbox{ and } \max_{1\le j \le 2} \kappa_{j}(x) < \frac 1R\mbox{ in }\  \mathbb R^2.
\end{equation}
Then, $c = F(1,1)$ and $f$ is an affine function, that is, $\partial\Omega$ must be a hyperplane.
\end{theorem}

\noindent
{\it Proof. } Here we have Lemma \ref{le:mean curvature inequality viscosity sense}. We consider $\Gamma^{*}$ defined by \eqref{middle parallel surface} as in Section \ref{Proof of 2nd theorem}. Then $\partial\Omega$, $\Gamma^{*}$, and $\Gamma$ are parallel hypersurfaces. Denote by $\kappa^{*}_{1}(Z), \kappa^{*}_{2}(Z)$ the principal curvatures of $\Gamma^{*}$ with respect to the upward unit normal vector $\nu^{*}(Z)$ to $\Gamma^{*}$
at $Z\in\Gamma^{*}$, and denote by $\hat{\kappa}_{1}(Y), \hat{\kappa}_{2}(Y)$ the principal curvatures of $\Gamma$ with respect to the upward unit normal vector at $Y = Z + \frac R2\nu^{*}(Z) \in \Gamma$. Also, here for the principal curvatures of $\partial\Omega$ we use the notation $\kappa_{1}(X), \kappa_{2}(X)$ instead of 
$\kappa_{1}(x), \kappa_{2}(x)$ with $(x,f(x))= X = Z- \frac R2\nu^{*}(Z) \in \partial\Omega$. These principal curvatures have the following relationship: 
$$
\kappa_{j}(X) = \frac {\kappa^{*}_{j}(Z)}{1+\frac R2\kappa^{*}_{j}(Z)}\ \mbox{ and }\ 
\hat{\kappa}_{j}(Y) = \frac {\kappa^{*}_{j}(Z)}{1-\frac R2\kappa^{*}_{j}(Z)}\ \mbox{ for each } j = 1, 2.
$$
Since $\displaystyle\max_{1\le j \le 2} \kappa_{j}(X) < \frac 1R$ and $\displaystyle1-R\kappa_{j}(X) = \frac 1{1+R\hat{\kappa}_{j}(Y)}$,
we see that
\begin{equation}
\label{curvature bounds from above and below}
-\frac 2R < \kappa^{*}_{j}(Z) < \frac 2R\ \mbox{ for each } j = 1, 2.
\end{equation}
On the other hand, by Lemma \ref{le:mean curvature inequality viscosity sense}, we have
$$
\sum_{j=1}^{2}\frac {\kappa^{*}_{j}(Z)}{1+\frac R2\kappa^{*}_{j}(Z)} \le 0 \le \sum_{j=1}^{2}\frac {\kappa^{*}_{j}(Z)}{1-\frac R2\kappa^{*}_{j}(Z)}.
$$
This gives
$$
\kappa^{*}_{1}+\kappa^{*}_{2} +R\kappa^{*}_{1}\kappa^{*}_{2} \le 0 \le \kappa^{*}_{1}+\kappa^{*}_{2} -R\kappa^{*}_{1}\kappa^{*}_{2},
$$
and hence
$$
\kappa^{*}_{1}\kappa^{*}_{2} \le 0 \mbox{ and } R\kappa^{*}_{1}\kappa^{*}_{2} \le \kappa^{*}_{1}+\kappa^{*}_{2} \le -R\kappa^{*}_{1}\kappa^{*}_{2}.
$$
Then, with the aid of \eqref{curvature bounds from above and below}, we conclude that 
$$
(\kappa^{*}_{1})^{2}+ (\kappa^{*}_{2})^{2} \le 2\cdot(-3) \kappa^{*}_{1}\kappa^{*}_{2}.
$$
Hence the Gauss map of $\Gamma^{*}$ is $(-3, 0)$-quasiconformal on $\mathbb R^{2}$ (see \cite[(16.88), p. 424]{GT}) and hence by \cite[Corollary 16.19, p. 429]{GT} $\Gamma^{*}$ must be a hyperplane. \qed



\begin{small}

\end{small}
\end{document}